\documentclass[a4paper,12pt]{article}
\usepackage[margin=1.93cm]{geometry}

\usepackage{authblk}
\usepackage{graphicx}       
\usepackage{float}          
\usepackage[dvipsnames]{xcolor}    
\usepackage{booktabs}
\usepackage{siunitx}
\usepackage{subfig}
\usepackage{xurl}
\usepackage{placeins}
\usepackage{amssymb}
\usepackage{amsthm}
\usepackage{mathtools}

\usepackage{algorithm}
\usepackage{algpseudocode}
\makeatletter
\algnewcommand{\AlgComment}[1]{\Statex \hskip\ALG@thistlm \textcolor{gray}{\(\triangleright\) #1}}
\makeatother

\usepackage{microtype}      

\usepackage{hyperref}       
\hypersetup{
  colorlinks=true,
  linkcolor=blue,
  citecolor=blue,
  urlcolor=blue
}
\usepackage[noabbrev,nameinlink,capitalize]{cleveref}
\crefname{equation}{Eq.}{Eqs.}
\crefname{appsection}{}{}
\Crefname{appsection}{}{}

\newcommand{\tildeR}{R}
\usepackage[utf8]{inputenc}
\usepackage{newunicodechar}
\newunicodechar{−}{-}

\date{\today}

\title{Exact and Evolutionary Algorithms for Sequential Multi-Objective Transmission Topology Planning}

\author[1]{Job Groeneveld}
\author[2]{Miguel Mu\~{n}oz}
\author[1]{Jan Viebahn\thanks{Corresponding author: jan.viebahn@tennet.eu}}
\author[3]{Alessandro Zocca}

\affil[1]{TenneT TSO B.V., Arnhem, The Netherlands}
\affil[2]{Artelys, Paris, France}
\affil[3]{Department of Mathematics, Vrije Universiteit Amsterdam, The Netherlands}

\makeatletter
\renewcommand{\title}[1]{}
\renewcommand{\author}[2][]{}
\newcommand{\address}[2][]{}
\newcommand{\ead}[1]{}
\newcommand{\corref}[1]{}
\newcommand{\cortext}[2][]{}
\let\arxiv@bibliographystyle\bibliographystyle
\renewcommand{\bibliographystyle}[1]{\arxiv@bibliographystyle{abbrv}}
\newif\ifarxiv@maketitled
\let\arxiv@abstract\abstract
\let\arxiv@endabstract\endabstract
\renewenvironment{abstract}{%
  \ifarxiv@maketitled\else
    \maketitle
    \global\arxiv@maketitledtrue
  \fi
  \arxiv@abstract
}{%
  \arxiv@endabstract
}
\newenvironment{frontmatter}{}{}
\newenvironment{keyword}{%
  \par\medskip\noindent\textbf{Keywords: }%
}{%
  \par\medskip
}
\newcommand{\sep}{;\ }
\makeatother

\begin{document}

\begin{frontmatter}

\title{Exact and Evolutionary Algorithms for Sequential Multi-Objective Transmission Topology Planning}

\author[tennet]{Job Groeneveld}
\author[artelys]{Miguel Mu\~{n}oz}
\author[tennet]{Jan Viebahn\texorpdfstring{\corref{cor1}}{}}
\ead{jan.viebahn@tennet.eu}
\author[vu]{Alessandro Zocca}

\cortext[cor1]{Corresponding author}

\address[tennet]{TenneT TSO B.V., Arnhem, The Netherlands}
\address[artelys]{Artelys, Paris, France}
\address[vu]{Department of Mathematics, Vrije Universiteit Amsterdam, The Netherlands}

\begin{abstract}
We study day-ahead transmission topology control for high-voltage grid operation under $N-1$ security constraints. The operational task is to select, over a 24-hour horizon, a sequence of substation topologies obtained via busbar-coupler switching to relieve line overloads while limiting switching effort and topological complexity. We formulate this task as a sequential multi-objective optimization problem with four objectives used in TSO decision making: worst-case $N-1$ line loading, maximum topological depth, number of topology changes, and time spent outside the reference topology. We propose an exact block algorithm that exploits the temporal structure of topology plans: consecutive hours with the same topology are represented as blocks, enabling enumeration of the complete Pareto front over the admissible set of topologies under fixed operational bounds on depth and switching. We also develop a tailored NSGA-III-based evolutionary heuristic and evaluate it against the exact front. Using real operational data from the Dutch high-voltage transmission grid operated by TenneT, the block algorithm computes the exact front for a highly congested day in under three minutes after topology-level load-flow preprocessing. The exact front reveals low-switching plans with no DC $N-1$ thermal overloads that the tested evolutionary search fails to find. The proposed method, therefore, provides both a practical day-ahead decision-support tool for transmission operators and a benchmark for heuristic and learning-based topology-control methods.
\end{abstract}

\begin{keyword}
Transmission topology control \sep $N-1$ security \sep Congestion management \sep
Optimal transmission switching \sep Busbar switching \sep Multi-objective optimization \sep
TSO decision support
\end{keyword}

\end{frontmatter}

\section{Introduction}
\label{Introduction}
Transmission system operators (TSO) routinely assess whether the day-ahead high-voltage grid schedule remains secure under $N-1$ contingency criteria~\cite{Hedman2009Contingency,Hedman2010UnitCommitment}. When forecasted injections lead to overloads, operators can deploy remedial actions to modify power flows before more expensive measures, such as redispatch, are required. One operationally important class of remedial actions is transmission topology control: changing the network configuration so that power is redistributed over the existing grid while relevant line-loading limits remain respected~\cite{HedmanOren2011Review,Numan2023TheReview}.

This paper focuses on substation-level topology control through busbar-coupler switching. Opening or closing busbar couplers changes the electrical grouping of busbars inside a substation and can be modeled as node splitting in the transmission network~\cite{Heidarifar2016}. This differs from line switching, which is the dominant abstraction in much of the optimal transmission switching literature~\cite{Fisher2008OTS,Hedman2009Contingency,Hedman2010UnitCommitment}. From a control-room perspective, busbar-coupler switching creates a multi-objective decision problem: a topology may reduce the worst $N-1$ line loading, but at the cost of additional switching actions, greater topological depth, or more time spent operating away from the reference topology.

The day-ahead version of this problem is inherently sequential. A topology chosen for one hour affects the number of subsequent switching actions, and the operator is not seeking a single mathematical optimum but a set of feasible trade-offs among security margin, switching effort, topological complexity, and persistence of non-reference configurations. Existing topology-control methods are often heuristic, single-objective, or applied independently at each time step~\cite{Granelli2006OptimalAlgorithms,Zhang2014OptimalReliability,Marot2020LearningControllers,Dorfer2022PowerAlphaZero}. As a result, they do not provide a complete view of the Pareto-optimal operational choices over the full planning horizon, nor do they provide an exact reference against which heuristic or learning-based methods can be evaluated.

Increasing demand and renewable-driven flow patterns push transmission grids to operate closer to their limits, which raises the value for TSOs of corrective actions that do not require redispatch. Topology control is attractive in this setting because it can relieve overloads without immediately requiring grid expansion or generation redispatch, provided that the resulting topology remains acceptable from an operational security perspective~\cite{HedmanOren2011Economic,Dorfer2022PowerAlphaZero}. The relevant question for a TSO is therefore not only whether a topology improves a loading metric, but whether a full day-ahead topology plan provides an acceptable trade-off between $N-1$ security, switching burden, and operational simplicity.

The key methodological gap is that sequential, multi-objective topology control has not been solved exactly at the level needed for day-ahead operational decision support. Existing exact formulations for optimal transmission switching primarily target single operating points or economic objectives~\cite{Fisher2008OTS,ONeill2010EconomicConcepts}. For a 24-hour horizon, direct enumeration over hourly topology choices is combinatorially prohibitive. This has encouraged the use of heuristics, dynamic programming approximations, evolutionary algorithms, and reinforcement-learning methods~\cite{Granelli2006OptimalAlgorithms,Zhang2014OptimalReliability,Marot2020LearningControllers,Dorfer2022PowerAlphaZero}. Such approaches can produce useful candidate actions, but without an exact Pareto front it remains difficult to know whether they miss operationally important trade-offs, especially low-switching strategies that remove DC $N-1$ thermal overloads.

We address this gap by exploiting the block structure of day-ahead topology plans. Instead of treating each hour independently, we group consecutive hours with the same topology into blocks and optimize over the resulting block configurations. This preserves the operational objectives exactly while reducing the effective search space for fixed bounds on switching and topological depth. We then use the exact Pareto front as both a decision-support object and a benchmark for a tailored multi-objective evolutionary algorithm (MOEA).

The main contributions of this paper are as follows:
\begin{enumerate}
    \item We formalize day-ahead substation-level transmission topology control as a sequential multi-objective optimization problem under $N-1$ security assessment, using four operational objectives reflecting real TSO decision criteria.
    \item We develop an exact algorithm, named \emph{block algorithm}, that exploits the temporal block structure of strategies to enumerate the complete Pareto front over the precomputed admissible topology set for fixed maximum depth and switch count.
    \item We design a tailored MOEA with structure-guided initialization and problem-specific variation operators, and evaluate it against the exact Pareto front using both $\mathrm{IGD}^+$ 
    and a novel multi-front coverage metric.
    \item We demonstrate, using real operational data from the Dutch high-voltage grid operated by TenneT, that exact enumeration is feasible and recovers operationally relevant plans missed by the tested evolutionary search, providing transmission operators with a practical and provably complete decision-support tool.
\end{enumerate}

The rest of the paper is organized as follows. \cref{sec:Literature_Review} reviews related work on optimal topology switching and optimization methods. \cref{sec:methodology} formulates the problem and solution algorithms. \cref{seq:experimental_setup} describes the case study and evaluation framework. \cref{sec:results} reports the exact and evolutionary results, and \cref{sec:conclusion} discusses operational implications, limitations, and future work.

\section{Related Literature}
\label{sec:Literature_Review}

Transmission topology control is commonly studied under the term optimal topology switching (OTS), which covers dynamic grid reconfiguration through changes in network topology. OTS formulations differ mainly in the power-flow model, the inclusion of flexibility options such as dynamic thermal rating, storage, or renewable curtailment, and the treatment of uncertainty~\cite{HedmanOren2011Review,Numan2023TheReview}. Most studies model line-switching decisions, whereas this paper addresses substation reconfiguration through busbar-coupler switching, or node splitting. Explicit node-breaker formulations have received less attention; a notable exception is~\cite{Heidarifar2016}. The present work uses deterministic input data, a DC load-flow model~\cite{10314781}, and busbar-coupler switching as the only control action.

Classical OTS research has primarily used mathematical programming formulations. The foundational work on OTS was introduced by~\cite{Fisher2008OTS}, who formulated the problem as a mixed-integer program (MIP) with binary variables representing the status of each transmission element and demonstrated substantial economic savings on the IEEE 118-bus system. 
Subsequent studies by Hedman et al.~\cite{Hedman2009Contingency,Hedman2010UnitCommitment} added $N-1$ contingency constraints and unit-commitment coupling, showing that topology changes can vary by hour while maintaining reliability standards. DC-based OTS is most often formulated as a mixed-integer linear program with an economic objective and thermal-limit constraints~\cite{HedmanOren2011Economic,Numan2023TheReview,Numan2020MobilizingIntegration,ONeill2010EconomicConcepts,Li2021ASwitching,Behnia2019IntegratedSwitching}.

Learning-based methods have also been applied to grid topology optimization, partly driven by RTE's \emph{Learning to Run a Power Network} challenge~\cite{Marot2020LearningControllers}; see~\cite{vanderSar2025} for an overview. Related work includes the AlphaZero approach of~\cite{Dorfer2022PowerAlphaZero}, which reduced redispatching volume but required substantial training and computation. These methods are promising for fast online decision support, but they do not generally provide an exact Pareto reference set for a given day-ahead operating condition.

Operational decision-support tools often combine optimization with domain-specific screening. TenneT has developed \textit{GridOptions}, a tool that combines dynamic programming with domain-specific heuristics~\cite{ViebahnGridOptions}. It first computes load flows for candidate topologies under load and generation forecasts, filters topologies by congestion-free duration, and then derives sequential strategies on a decision graph using random-weight sampling and Dijkstra's algorithm.

The four operational objectives used in the present study (see~\cref{sec:problem_objectives}) were originally identified in the context of GridOptions as key performance metrics for evaluating topological strategies in TenneT's day-ahead congestion management workflow~\cite{ViebahnGridOptions}. However, a head-to-head numerical comparison with GridOptions is outside the scope of this work.

Evolutionary algorithms have also been applied to OTS. \cite{Granelli2006OptimalAlgorithms} used a genetic algorithm for single-objective OTS, while \cite{Zhang2014OptimalReliability} used a multi-objective evolutionary algorithm to approximate trade-offs between generation cost and loss-of-load probability.
More broadly, MOEAs such as NSGA-II are widely used for combinatorial optimization in power systems~\cite{Ma2023AApplications}. Because the present problem has four competing objectives, it falls in the many-objective setting, for which NSGA-III~\cite{Deb2014AnConstraints} is adopted as the selection mechanism (see~\cref{sec:MOEA_algo}).

Taken together, this literature shows that topology control is operationally relevant but still lacks an exact sequential, multi-objective treatment at the level needed for day-ahead TSO decision support. To the best of our knowledge, no exact algorithm has been proposed for sequential, multi-objective OTS. Existing exact approaches primarily target single economic objectives~\cite{Fisher2008OTS,ONeill2010EconomicConcepts}, whereas heuristic and learning-based methods improve scalability without global optimality guarantees. Two gaps are therefore central here: operational topology control has not been solved exactly for simultaneous security, switching, depth, and reference-topology objectives, and the temporal block structure of day-ahead strategies has not been exploited to enumerate all Pareto-optimal choices. The block algorithm addresses both gaps.

\section{Methodology}
\label{sec:methodology}
This section introduces the mathematical formulation of the congestion management problem and the two proposed solution methods. \cref{sec:sets_parameters} defines the fundamental terminology and mathematical notation. Based on these definitions, \cref{sec:problem_objectives} formulates the optimization objectives. Finally, \cref{sec:block_algo,sec:MOEA_algo} describe the proposed block algorithm and MOEA, respectively. We present here the algorithmic structure needed to understand the results. \cref{app:blockalgorithmcomplexity} gives the complexity proof for the block algorithm, and \cref{app:moea} gives the MOEA representation, hyperparameters, and evaluation metrics.

\subsection{Preliminaries}\label{sec:sets_parameters}
We begin by introducing the terminology and notation used in the rest of the paper. 
\begin{itemize}
    \item We consider a planning horizon consisting of $T_{\max}$ unit periods. We identify these unit periods with their starting times in $T :=\{0, 1, 2, \dots, T_{\max}-1\}$, which we refer to as \textbf{time steps}\footnote{In our implementation, we consider a one-day time horizon, with each unit period corresponding to one hour, hence featuring $|T|=T_{\max}= 24$ time steps. However, the framework presented here is general and can readily be adapted to different time resolutions and time horizons.}.
    
    \item $G_t$ denotes the collection of available network topologies at time step $t \in T$ (and remains so for the whole following unit period, until the next time step). Each available topology $g \in G_t$ is described by a heterogeneous graph $g = (BB(g), BR(g), PO(g), E(g))$, with node set $BB(g) \cup BR(g) \cup PO(g)$ and edge set $E(g)$, where $BB(g)$ denotes the set of busbars, $BR(g)$ the set of transmission branches (i.e., lines and transformers), and $PO(g)$ includes the set of power injections (i.e., load and generation). The edges in $E(g)$ represent how the nodes are connected.
    
    \item The \textbf{reference topology} $\tilde{g}$ is the graph corresponding to the network configuration with all busbar couplers closed, which is assumed to be always available, i.e., $\tilde{g} \in G_t$ for every $t \in T$.

    \item For each time step $t \in T$, we consider a set of forecasted power injections. Although stochastic in nature, these are treated as fixed and deterministic parameters during optimization. If topology $g \in G_t$ is implemented, branch flows are derived using these injections and the DC power flow approximation~\cite{10314781}. For each branch $b \in BR(g)$, we denote by $u_b(g, t)$ the \textit{capacity utilization}, calculated as the ratio of the branch flow to its thermal limit under the worst-case $N-1$ contingency scenario. Capacity utilization is expressed as a percentage and may exceed 100\% in the event of a line overload.
    
    \item For each time step $t \in T$ and each topology $g \in G_t$, the worst-case $N-1$ line loading is given by $LF_1(g,t) \coloneq \max_{b \in BR(g)} u_b(g,t)$.
    
    \item For each topology $g$, the \textit{topological depth} $d(g) \coloneq |BB(g)| - |BB(\tilde{g})|$ counts the number of additional busbars created by opening busbar couplers relative to the reference topology $\tilde{g}$.
    
    \item A \textit{strategy} $s = (g_t)_{t \in T} \in S$ is an ordered 
    sequence of topologies, specifying which topology $g_t \in G_t$ is 
    deployed at each time step $t \in T$. The set of all possible 
    strategies is obtained as the Cartesian product of the available 
    topologies across all time steps, that is, 
    $S \coloneq \bigtimes_{t \in T} G_t$.
\end{itemize}

For each topology, time step, and monitored contingency, the DC load-flow calculation solves for the active power flows induced by the forecasted injections on the corresponding network topology. Branch loading is then normalized by the relevant thermal limit, and $LF_1(g,t)$ is obtained by taking the maximum utilization over monitored branches and single-element contingencies. Thus, an $LF_1$ value below $100\%$ means that no DC thermal overload is observed for the considered $N-1$ contingency set; it does not by itself certify AC feasibility, voltage-security margins, dynamic stability, or protection constraints.

\subsection{Problem objectives}
\label{sec:problem_objectives}
TenneT identified four key performance metrics for evaluating topological strategies~\cite{ViebahnGridOptions}, all of which we aim to minimize. Specifically, for any strategy $s = (g_t)_{t \in T} \in S$, we consider its:

\begin{enumerate}
    \item \textbf{Worst $N-1$ line capacity utilization across the entire time horizon}: 
    The first objective seeks to minimize the worst-case network loading across all time steps and all single-line contingencies, that is,
    $
    LF_1(s) := \max_{t \in T} LF_1(g_t,t);
    $
    \item \textbf{Topological depth across the entire time horizon}: 
    Node splitting increases the network complexity~\cite{Zhou2021}, so from the operator perspective, it is convenient to minimize the number of open busbar couplers, motivating the following choice for the second objective:
    $
    d(s) := \max_{t \in T}  d(g_t);
    $
    \item \textbf{Number of topology switches over the entire time horizon}: 
    Switching actions cause mechanical wear on switching assets and carry a small but non-negligible risk of equipment failure or operational error, potentially reducing the system reliability~\cite{Pranji2024}. For this reason, we introduce a third objective to track the number of topology switches required by a strategy as
    $
    w(s) \coloneq \sum_{t \in T \setminus \{0\}} \mathbf{1}_{\{g_t \neq g_{t-1}\}};
    $
    \item \textbf{Utilization of non-reference topologies}: 
    For network operators, it is often desirable to minimize the time spent on network configurations that differ from the reference topology $\tilde{g}$, in which all busbars are closed. It is thus natural to formulate a fourth objective as 
    $
    z(s) \coloneq \sum_{t \in T} \mathbf{1}_{\{g_t \neq \tilde{g}\}},
    $
    to count across the time horizon in how many time steps the reference topology $\tilde{g}$ is not used.
\end{enumerate}

Since all four objectives are minimized simultaneously, we compare strategies using Pareto dominance. 
For two strategies $s,s' \in S$, define
\[
s' \preceq s
\quad \Longleftrightarrow \quad
LF_1(s') \le LF_1(s),\quad d(s') \le d(s),\quad
w(s') \le w(s),\quad z(s') \le z(s),
\]
and $s' \prec s$ if and only if $s' \preceq s$ and at least one of the four inequalities is strict. For any feasible strategy set $A \subseteq S$, its Pareto front is
$
PF(A) := \{\, s \in A : \nexists\, s' \in A \text{ such that } s' \prec s \,\}.
$
In the numerical results, we distinguish between Pareto-optimal strategies and unique objective-space points, since several structurally distinct topology plans may yield identical objective values.

\subsection{Block algorithm} 
\label{sec:block_algo}
This subsection introduces the block algorithm for the multi-objective problem in \cref{sec:problem_objectives}. We define the block representation, show how feasible strategies are recovered, and give the algorithmic logic; the complexity analysis is deferred to \cref{app:blockalgorithmcomplexity}.
The key idea behind the block algorithm is that, from an operational perspective, only topology changes matter, not the exact time steps at which identical topologies persist. By grouping consecutive time steps with identical topologies into blocks, we drastically reduce the effective decision space while preserving exact optimality.

The block algorithm represents strategies by partitioning the time horizon into contiguous, non-overlapping intervals called blocks, during which the network topology remains constant. Formally, we define a \textbf{block} $B$ to be a set of consecutive time steps,
$
    B \coloneq \{ t_s, t_s +1, \dots, t_e \},
$
where $t_s,t_e \in T$ and $t_s \leq t_e$. Its length is
$
|B|:=t_e-t_s+1.
$
Let $\mathcal{B}$ be the collection of all such blocks. Blocks $B,B' \in \mathcal{B}$ are consecutive if $t_s(B') = t_e(B) + 1$.

A \textbf{block configuration} is a collection $M=\{B_1,\dots,B_{|M|}\}$ of non-overlapping blocks that covers the horizon, i.e., $B_i \cap B_j=\emptyset$ for $i\ne j$ and $\bigcup_{i=1}^{|M|}B_i=T$. Let $\mathcal{M}(l)$ be the set of block configurations with $l+1$ blocks:
\begin{equation}
\label{eq:block_conf_sets}
    \mathcal{M}(l) \coloneq \{\, M \subseteq \mathcal{B} ~:~ M \text{ is a block configuration and} \, |M| = l + 1 \,\}\ .
\end{equation}
Each configuration in $\mathcal{M}(l)$ partitions $T$ into $l+1$ consecutive blocks.


The set of topologies always available within a block $B$ is 
$
G_B \coloneq \bigcap_{t \in B} G_t.
$
For any $G \subseteq G_B$, the best load-flow value attainable by a constant topology over block $B$ is
    $
    LF_1(B, G) \coloneq \min_{g \in G} \max_{t \in B} LF_1(g, t),
    $
with the convention that $LF_1(B,G)=\infty$ if $G$ is empty. If the second argument is omitted, the minimum is taken over all topologies available in that block, i.e.,
    $
    LF_1(B) := LF_1(B, G_B).
    $
A natural subset of topologies available in block $B$ to which we might want to restrict is the subset $G_B^{(\leq d)}$ of topologies available in block $B$ that have a topological depth at most $d$, i.e.,
    $
    G_B^{(\leq d)} \coloneq \{ g \in G_B ~:~ d(g) \leq d \,\} \subseteq G_B.
    $
This is particularly relevant for objective 2 in \cref{sec:problem_objectives}. The best achievable $LF_1$ metric for block $B$, restricted to topologies with depth at most $d$, is denoted by
$
LF_1(B, d) \coloneq LF_1 (B, G_B^{(\leq d)}).
$
The notation for best achievable $LF_1$ metrics extends from blocks to block configurations: 
$
    LF_1(M) \coloneq \max_{B \in M} LF_1(B)
$
is the best achievable $LF_1$ for a block configuration $M$, and
\begin{equation}\label{eq:LF1Md}
    LF_1(M, d) \coloneq \max_{B \in M} LF_1(B, d)
\end{equation}
is the best achievable $LF_1$ for a block configuration $M$ restricted to topologies of depth at most $d$. 

A strategy $s \in S$ can now be described by assigning one topology to each block. For a strategy $s=(g_t)_{t\in T}$ and block $B$, write $s_{|B}=g$ if $g_t=g \in G_B$ for all $t\in B$. The strategies whose coarsest block representation is compatible with a block configuration $M$ are
\[
    S(M) \coloneq \{ s \in S ~:~ s_{|B} = g \text{ for some } g \in G_B, \, \forall \, B \in M \}.
\]
Consecutive blocks are assumed to have distinct topologies, so the number of switching time steps for any strategy in $S(M)$ is equal to
$
w(M) \coloneq |M| -1,
$
by construction, since it is determined by the number of blocks in $M$. 

For objective~4, we track which blocks use the reference topology. For a block configuration $M$, let $\mathcal{A}(M)$ be the collection of non-consecutive reference-block assignments:
\begin{equation} \label{eq:block2}
\mathcal{A}(M) \coloneq \Bigl\{ \tildeR \subseteq M ~:~ \forall\, B,B' \in \tildeR : t_s(B) \neq t_e(B')+1, \,t_s(B') \neq t_e(B)+1 \Bigr\}.
\end{equation}
The restriction to non-consecutive reference blocks avoids representing two adjacent reference blocks that should instead be merged. For $\tildeR\in\mathcal{A}(M)$, define the subset $S(M, \tildeR)$ of strategies for which the utilization of the reference topology is explicit in all the blocks in $\tildeR$, that is $S(M, \tildeR) \coloneq \left \{ s \in S(M) ~:~ s_{|B} = \tilde{g} \, \forall \, B \in \tildeR\ \text{and } s_{|B} \neq \tilde{g} \, \forall \, B \in M \setminus \tildeR \right \}$. Note that, given that $\tilde{g} \in G_t$ for all $t \in T$ (as noted in~\cref{sec:sets_parameters}), it follows by definition that $\tilde{g} \in G_B$ for every block $B \in \mathcal{B}$.

The number $z(s)$ of time steps for which the reference topology is not used for any strategy $s \in S(M,\tildeR)$ is entirely determined by the blocks in $\tildeR$, which justifies the following notation
$
    z(\tildeR) := T_{\max} - \sum_{B \in \tildeR} |B|.
$
Given $(d, M, \tildeR)$, the best achievable $LF_1$ value is
\[
    LF_1(d, M, \tildeR) \coloneq \max \left\{ \max_{B \in M \setminus \tildeR} LF_1(B,G_B^{(\leq d)}), \,  \max _{B \in \tildeR} LF_1(B,\{\tilde{g}\}) \right\}.
\]
For each evaluated configuration $(d,M,\tildeR)$, the associated feasible strategies are obtained by filtering the available non-reference topologies in every non-reference block by depth and $LF_1$ threshold, and then taking the admissible assignments from the corresponding Cartesian product:
\begin{equation}
\label{eq:cartasian_product}
\mathcal{S}(d, M, \tildeR)
 \coloneq \left\{\big( (g_B)_{B \in M \setminus \tildeR},\,  (\tilde{g}_B)_{B \in \tildeR} \big) \,:\, d(g_B) \leq d, \,  \max_{t \in B} LF_1(g_B, t) \le LF_1(d, M, \tildeR) \ \forall B \in M \setminus \tildeR \right\}.
\end{equation}

Each strategy $s \in \mathcal{S}(d, M, \tildeR)$ specifies a feasible combination of topologies having the objective values $(d, w(M), z(\tildeR), LF_1(d, M, \tildeR))$. 

The block algorithm in \cref{alg:block} takes a maximum depth $d^{\max}$ and maximum number of switching time steps $s^{\max}$. It enumerates every depth level $d$, every block configuration $M\in\mathcal{M}(l)$ for $0\le l\le s^{\max}$, and every reference-block assignment $\tildeR\in\mathcal{A}(M)$. The resulting objective tuples are collected in $C$ and filtered for nondominance.

\begin{algorithm}[H]
\begin{algorithmic}[1]
    \Require Maximum topological depth $d^{\max}$, maximum number of switching time steps $s^{\max}$
    \AlgComment{Initialize set of objective values}
    \State $C \gets \emptyset$
    \For {$d = 0  \to d^{\max}$} 
        \For {$l = 0 \to s^{\max}$} 
            \AlgComment{Loop through every block configuration $M$}
            \For {$M \in \mathcal{M}(l)$}   
                \AlgComment{Loop through every possible reference topology assignment $\tildeR$}         
                \For {$\tildeR \in \mathcal{A}(M)$}
                    \AlgComment{Append objective values to set $C$:}
                    \State $C \gets C \cup (d, M, \tildeR, w(M), z(\tildeR), LF_1(d, M, \tildeR)) $
                \EndFor
            \EndFor
        \EndFor
    \EndFor
    \AlgComment{Filter nondominated configurations using Algorithm 2}
    \State $C^* \gets \text{FilterNonDominated}(C)$
    \State \Return $C^*$ 
\end{algorithmic}
\caption{BlockAlgorithm}
\label{alg:block}
\end{algorithm}

Each element of $C^*$ corresponds to at least one Pareto-optimal feasible strategy with respect to the four objectives. \cref{alg:nondominated_sorting} exploits the fact that all objectives except $LF_1$ take finite integer values: nondominated solutions are identified by comparing $LF_1$ values against neighboring entries in the discrete objective dimensions. This approach identifies nondominated solutions more efficiently than conventional sorting methods.

The rounding in \cref{alg:nondominated_sorting} avoids distinguishing solutions whose $LF_1$ values differ by less than 0.1, so nondominated sorting reflects operationally meaningful differences rather than numerical noise. The algorithm is exact for the stated depth and switching bounds; the combinatorial count of block configurations, the reduction obtained by excluding consecutive reference blocks, and the resulting polynomial scaling for fixed operational bounds are derived in~\cref{app:blockalgorithmcomplexity}.

\begin{algorithm}[H]
\begin{algorithmic}[1]
\Require Set of objective values $C$, maximum depth $d^{\max}$, maximum number of switching time steps $s^{\max}$, and set of time steps $T$.
\AlgComment{Initialize the tensor $O$ by setting all elements to infinity. A padding of $-1$ is added so that neighbor lookups are well-defined and do not require specific boundary cases}
\State $O_{d, l, n} \gets \infty \quad 
\forall \, d \in \{-1, 0, \dots, d^{\max}\}, \quad
\forall \, l \in \{-1, 0, \dots, s^{\max}\}, \quad 
\forall \, n \in  \{-1 , 0, \dots, T_{\max}\}$
\AlgComment{Store the best loading for each discrete objective triple $(d, w(M), z(\tildeR))$ from the set $C$}
\State $O_{d, l, n} \gets \min \{LF_1(d, M, \tildeR): {(d, M, \tildeR, w(M), z(\tildeR), LF_1(d, M, \tildeR)) \in C,}\ w(M)=l,\ z(\tildeR)=n\}$
\AlgComment{Initialize $R^{nds}$ as an empty collection of nondominated solutions}
\State $R^{nds} \gets \emptyset$
\AlgComment{Iterate over all elements in the tensor $O$}
\For{$d = 0 \to d^{\max}$}
    \For{$l = 0 \to s^{\max}$}
        \For{$n = 0 \to T_{\max}$}
            \AlgComment{Determine the minimum $LF_1$ value from neighboring elements}
            \State $LF_1^{nb} \gets \min\{O_{d-1,l,n},\; O_{d,l-1,n},\; O_{d,l,n-1}\}$
            \AlgComment{Compare the current $LF_1$ value after rounding to the first decimal place}
            \If{$\mathrm{round}(O_{d,l,n}) < \mathrm{round}(LF_1^{nb})$}
                \AlgComment{Store the current $LF_1$ value if it is lower than its neighboring values}
                \State $R^{nds} \gets R^{nds} \cup \{ (d, l, n, O_{d,l,n}) \}$
            \Else
                \AlgComment{Otherwise, update the current value to the best neighboring $LF_1$ value}
                \State $O_{d,l,n} \gets LF_1^{nb}$
            \EndIf
        \EndFor
    \EndFor
\EndFor
\State \Return $R^{nds}$
\end{algorithmic}
\caption{FilterNonDominated}
\label{alg:nondominated_sorting}
\end{algorithm}

\subsection{Multi-objective evolutionary algorithm (MOEA) framework}
\label{sec:MOEA_algo}

The multi-objective evolutionary algorithm (MOEA) framework is included for two purposes. First, it provides a representative heuristic topology-control method of the type that could be used when exact enumeration is unavailable. Second, because the block algorithm returns the exact Pareto front for the present instance, the MOEA can be evaluated against a known reference rather than only through indicators (such as hypervolume) that score the MOEA against itself. This comparison is important for TSO decision support: a heuristic that returns plausible topology plans may still miss low-switching plans that remove DC $N-1$ thermal overloads and are operationally valuable.

Each MOEA individual represents a complete day-ahead topology strategy $s=(g_t)_{t\in T}$. The gene at time step $t$ is one admissible topology $g_t\in G_t$, so feasibility with respect to topology availability is enforced directly by the encoding. The objective functions are exactly those in \cref{sec:problem_objectives}; no scalar weighting of $N-1$ loading, topological depth, switching effort, or non-reference operation is used. This is consistent with the operational interpretation of the output as a set of trade-off plans rather than a single optimized schedule.

The evolutionary design is adapted to the structure of the power-system problem. The initial population is stratified by switch count and topological depth because uniform sampling from the available topologies would strongly overrepresent high-depth strategies and would typically create unrealistic hour-to-hour switching. Variation is performed using $k$-point crossover, which preserves contiguous temporal segments of a topology plan, and random-reset mutation, which replaces selected hourly topologies by other admissible topologies at the same time step. Selection is performed using NSGA-III~\cite{Deb2014AnConstraints}, which is well-suited to this four-objective setting, as it combines nondominated sorting with reference-direction-based diversity preservation.

\Cref{app:moea} gives details on chromosome representation, structure-guided initialization, crossover, mutation, selection, hyperparameters, random seeds, and performance metrics.

\section{Experimental Setup}
\label{seq:experimental_setup}
This section outlines the dataset and protocol used to evaluate the block algorithm and to compare the MOEA with the exact reference front. We first introduce the operational dataset and its construction in \cref{sec:dataset}. The MOEA is run for 15 independent random seeds under small, medium, and large population settings, with mutation probabilities between 0.05 and 0.20 and a 30-minute computational budget. The complete hyperparameter table and the formal definitions of the evaluation metrics are given in \cref{app:hyperparameters,app:evaluation measures}.

\subsection{Dataset} 
\label{sec:dataset}
To obtain a transparent and computationally efficient evaluation, we decouple the power grid calculations from the optimization process by constructing an offline dataset. For each time step, this dataset contains the set of available topologies together with their associated congestion-related metrics, defining the complete search space explored by both the block and evolutionary algorithms.

The dataset is derived from real-world operational data from TenneT and represents a portion of the Dutch high-voltage transmission grid. All load-flow calculations are based on a single operational day selected for its high network congestion. As shown in~\cref{tab:timesteps}, the $LF_1$ values of the reference topology exceed $100\%$ for the majority of time steps, with peaks above $125\%$ during the evening hours ($t = 18$--$20$), indicating severe overloading under $N-1$ contingency conditions. This day, therefore, represents stressed operating conditions in which topology-based congestion management is most likely to add operational value over the reference dispatch.

\begin{table}[!ht]
\centering
\begin{tabular}{l|cccccccc}
\textbf{Time step $t$} & 0 & 1 & 2 & 3 & 4 & 5 & 6 & 7 \\
\midrule
$|G_t|$ & 39180 & 39180 & 39180 & 39180 & 39180 & 39180 & 39180 & 19950 \\
${LF_1}(\tilde{g})$ & 112.0 & 100.0 & 101.4 & 104.2 & 106.0 & 113.6 & 144.4 & 119.8 \\
\midrule
\textbf{Time step $t$} & 8 & 9 & 10 & 11 & 12 & 13 & 14 & 15 \\
\midrule
$|G_t|$ & 19590 & 19590 & 19590 & 19590 & 19590 & 19590 & 19590 & 19590 \\
${LF_1}(\tilde{g})$ & 126.2 & 120.0 & 115.6 & 118.0 & 118.8 & 113.3 & 101.4 & 103.2 \\
\midrule
\textbf{Time step $t$} & 16 & 17 & 18 & 19 & 20 & 21 & 22 & 23 \\
\midrule
$|G_t|$ & 19590 & 19590 & 19590 & 19590 & 19590 & 39180 & 39180 & 39180 \\
${LF_1}(\tilde{g})$ & 102.4 & 111.6 & 128.9 & 142.1 & 136.7 & 131.1 & 123.6 & 113.3 \\
\end{tabular}
\caption{Number of available topologies $|G_t|$ and load flow value ${LF_1}(\tilde{g})$ of the reference topology for each $t \in T$.}
\label{tab:timesteps}
\end{table}

To enable efficient large-scale evaluation, active power flows are computed using the DC approximation~\cite{10314781}, a standard linearization that assumes small angle differences, flat voltage profiles, negligible line resistance, and lossless transmission. While less accurate than full AC models, the DC approximation is widely used in transmission-level studies and provides a favorable trade-off between computational efficiency and modeling fidelity. In particular, it allows the repeated evaluation of hundreds of thousands of candidate topologies, which would be computationally infeasible using full AC power flow within the block algorithm or an evolutionary optimization loop.

Using DC power flow calculations, the capacity utilization $LF_1(g,t)$ and topological depth $d(g)$ are computed for all available topologies at each time step, as defined in \cref{sec:sets_parameters}. The dataset contains $252{,}094$ unique topologies in total across the 24 time steps. \cref{tab:timesteps} reports the count $|G_t|$ of available topologies per time step, together with the $LF_1$ value of the reference topology (${LF_1}(\tilde{g})$).

The distribution of topological depths across the dataset is summarized in \cref{tab:values_counts}. The dataset is highly imbalanced, with the majority of available topologies concentrated at depth 3.

\begin{table}[!ht]
\centering
\begin{tabular}{c|c|c|c|c}
\textbf{Topological depth } & \, 0 & 1 & 2 & 3 \\
\hline
\textbf{\# of topologies } & \, 1 & 489 & 20777 & 230827 \\
\end{tabular}
\caption{Number of unique available topologies at each depth $d=0,1,2,3$ in the dataset.}
\label{tab:values_counts}
\end{table}
\FloatBarrier

\section{Results}
\label{sec:results}
This section reports the results from an operational decision-support perspective. We first analyze the exact Pareto front produced by the block algorithm: the complete menu of nondominated day-ahead topology plans over the admissible topology set and imposed depth and switching bounds. We then compare this reference set with the MOEA output using $\mathrm{IGD}^+$ and two dominance-front coverage indicators, $I_k(P)$ and $\hat{I}_k(P)$. Formal metric definitions and the MOEA configuration table are provided in \cref{app:hyperparameters,tab:algorithm_configurations}.

\begin{figure}[!t]
    \centering
    \includegraphics[width=0.8\linewidth]{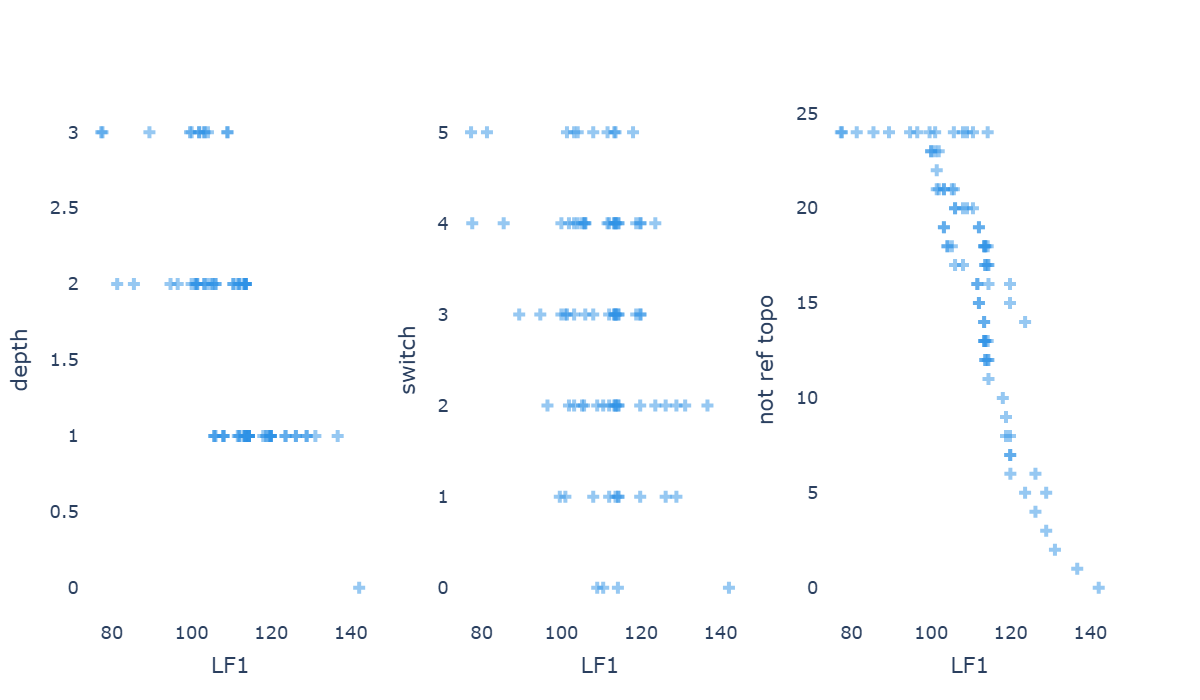}
    \caption{Pareto-optimal solutions for the block algorithm, plotted against topological depth (left), number of switching time steps (center), and time steps outside the reference topology (right).}
    \label{fig:distribution_block_only}
\end{figure}

\cref{fig:distribution_block_only} shows the objective values of the strategies identified by the block algorithm. Three key observations can be made:
\begin{itemize}
    \item A strategy with a maximum $LF_1$ value of approximately $80\%$ can be achieved with a maximum topological depth of $3$, requiring five switches and without relying on the reference topology.
    \item For the selected day, strategies achieving $LF_1 < 100\%$ are only attainable when the reference topology is not used.
    \item In total, the block algorithm identifies $163{,}807{,}652$ Pareto-optimal strategies, of which $83$ correspond to unique points in objective space. This indicates that identical Pareto-optimal trade-offs can be realized by structurally distinct strategies.
\end{itemize}
From an operational perspective, severe DC $N-1$ thermal overloads can be eliminated only by sustained non-reference operation at higher topological depth.

For a TSO, the exact front separates feasibility from preference: operators can first identify plans that remove the considered DC $N-1$ thermal overloads and then choose among them by switching effort, topological depth, and time away from the reference topology. The MOEA comparison below is therefore interpreted as a benchmark against this complete decision menu, rather than as an independent replacement for it.

\begin{figure}[!t]
    \centering
    \includegraphics[width=0.8\linewidth]{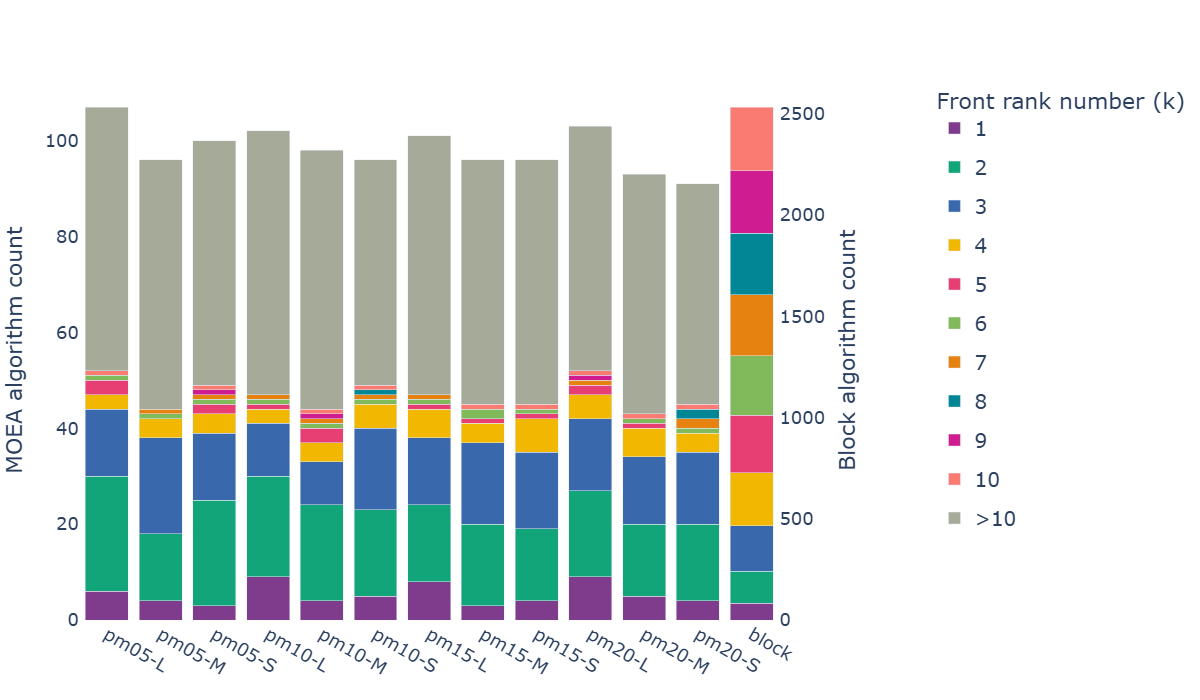}
   \caption{Distinct objective-space points per dominance front rank $k$, where lower rank fronts correspond to better trade-offs. For each MOEA configuration, the bars represent the count of recovered solutions per front, defined by $I_k(P) $ (left y-axis, see~\cref{eq:evaluation_metrics}). For the block algorithm, the rightmost bar reports the size of each reference front $|F_k|$ (right y-axis).}
    \label{fig:stacked_bar}
\end{figure}

\cref{fig:stacked_bar} compares the final combined Pareto front from each MOEA configuration with the block-algorithm reference fronts. The left y-axis reports the number of MOEA solutions assigned to each reference front, while the right y-axis reports the reference-front sizes; gray segments denote solutions outside the first ten fronts.

The block algorithm identifies substantially more unique solutions than any MOEA configuration. Among the MOEA variants, larger population sizes consistently yield more solutions within the first 10 dominance fronts. In contrast, variations in mutation probability do not appear to significantly influence the distribution of recovered fronts.

\begin{figure}[!t]
    \centering
    \includegraphics[width=0.8\linewidth]{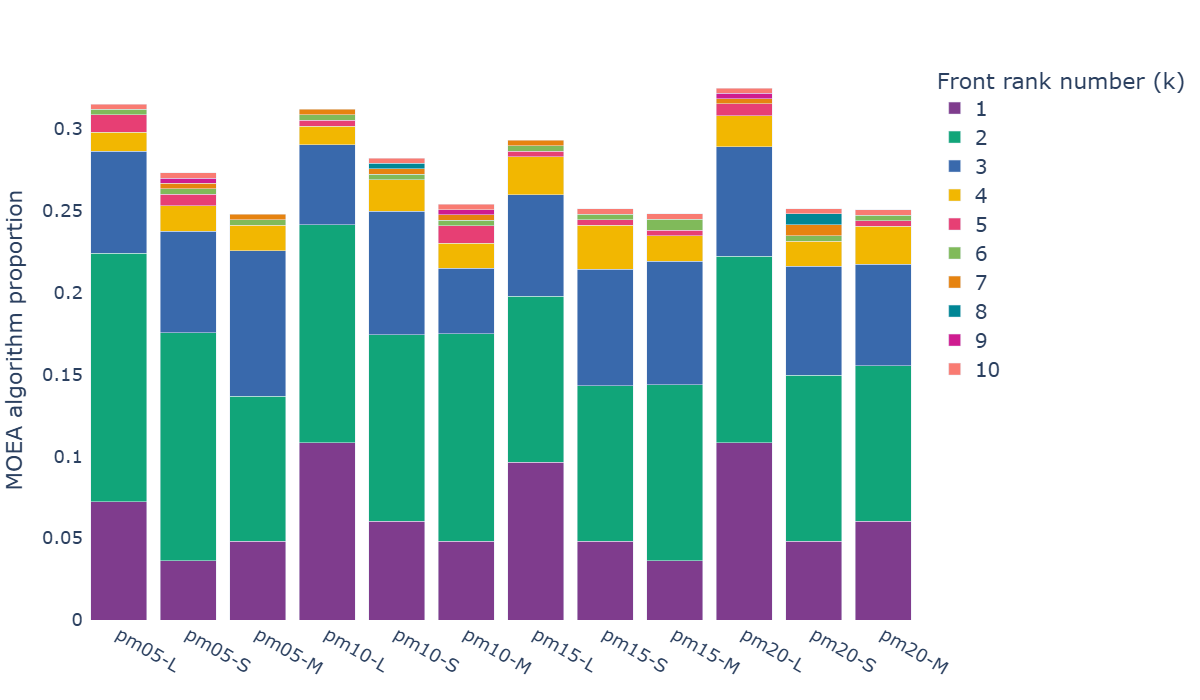}
   \caption{Relative coverage of reference dominance fronts across MOEA configurations for each reference front $F_k$. Segments represent the proportion of unique solutions recovered from each front $F_k$, see metric $\hat{I}_k(P)$ in~\cref{eq:evaluation_metrics}.}
    \label{fig:stacked_bar_proportion}
\end{figure}

\cref{fig:stacked_bar_proportion} shows that, across configurations, the MOEA recovers about $6\%$ of the first Pareto front and roughly $3$--$4\%$ cumulatively over the first ten dominance fronts. Coverage is highest on the second dominance front, at approximately $15\%$.
This limited coverage is partly attributable to the restricted population sizes relative to the size of the reference fronts.

\begin{figure}[!t]
    \centering
    \includegraphics[width=0.8\linewidth]{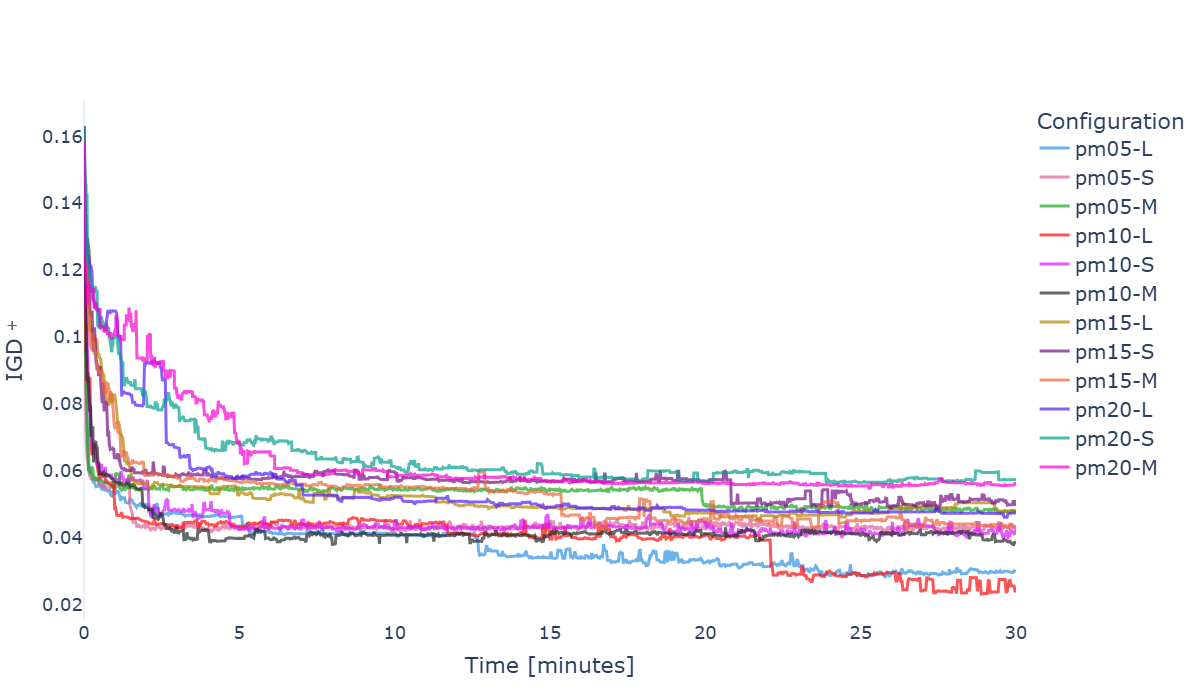}
    \caption{$\mathrm{IGD}^+$ trajectories for all MOEA configurations over runtime. The x-axis denotes runtime in minutes, while the y-axis shows the normalized $\mathrm{IGD}^+$ value.}
    \label{fig:igd_trajectories}
\end{figure}

\cref{fig:igd_trajectories} presents the $\mathrm{IGD}^+$ trajectories for all MOEA configurations. Convergence is fastest during the first $5$ minutes and largely saturates after approximately $25$ minutes. Configuration pm10-L attains the lowest $\mathrm{IGD}^+$ values and is therefore selected for further analysis.

\begin{figure}[!t]
    \centering
    \includegraphics[width=0.8\linewidth]{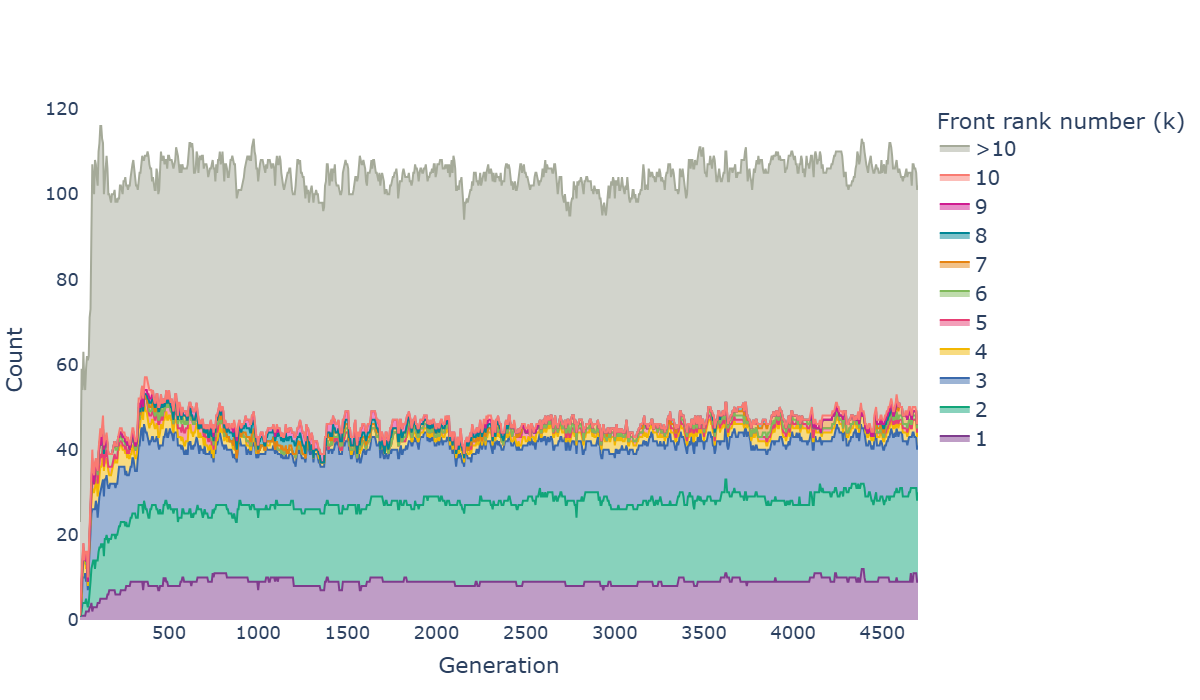}
    \caption{Evolution of dominance-front distribution for pm10-L. Lower rank fronts correspond to better trade-offs. The x-axis represents the generation number, and the y-axis indicates the number of solutions. Each colored area corresponds to one of the first ten dominance fronts, while the gray area represents solutions outside the top ten.}
    \label{fig:front_distribution_pm10-L}
\end{figure}

\cref{fig:front_distribution_pm10-L} shows that the dominance-front distribution for pm10-L stabilizes after approximately $1{,}500$ generations, with roughly half of the population in the first ten fronts. The $\mathrm{IGD}^+$ drop around minute 23 is not visible as a change in front counts because a newly discovered point can improve distance to the reference set without changing the number of first-front points.

\begin{figure}[!ht]
\label{fig_final_objectives}
    \centering
    \includegraphics[width=0.8\linewidth]{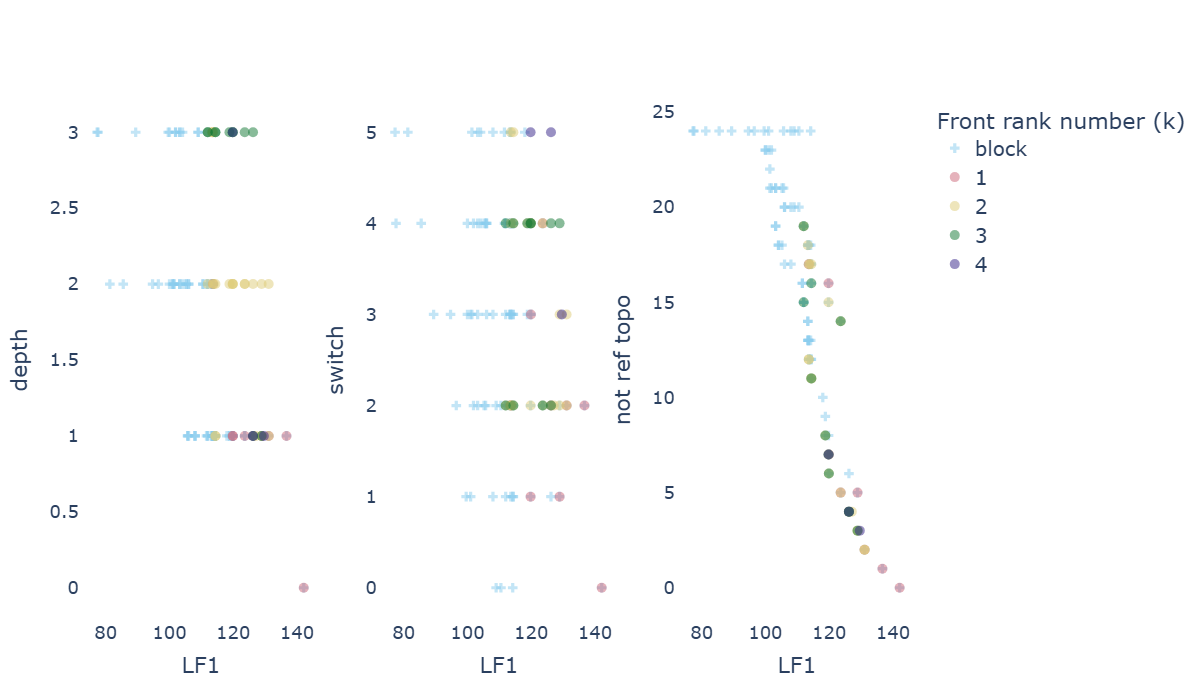}
    \caption{Objective-space distribution of final solutions for pm10-L, plotted against topological depth (left), number of switches (center), and time steps outside the reference topology (right).}
    \label{fig:front_distribution_conf010}
\end{figure}

\cref{fig:front_distribution_conf010} shows the final distribution of solutions across the first four dominance fronts for pm10-L. While the block algorithm identifies strategies with $LF_1<100\%$ at depth 2 or with a single switch, the MOEA does not find any solutions with $LF_1 < 100\%$.

\begin{table}[!h]
\centering
{\centering\textsc{Per front}\par}
\smallskip

\begin{tabular}{r rr rr}
\toprule
& \multicolumn{2}{c}{\textbf{MOEA} (pm10-L)} 
& \multicolumn{2}{c}{\textbf{Block algorithm}} \\
\cmidrule(lr){2-3} \cmidrule(lr){4-5}
\textbf{Front} $k$ 
& \textbf{Points} & \textbf{Strategies} 
& \textbf{Points} & \textbf{Strategies} \\
\midrule
 1 &  9 &  23 &    83 & $1.638 \times 10^{8\phantom{0}}$ \\
 2 & 21 &  54 &   158 & $1.380 \times 10^{13}$ \\
 3 & 13 &  30 &   225 & $5.394 \times 10^{15}$ \\
 4 &  5 &   7 &   261 & $4.714 \times 10^{16}$ \\
 5 &  2 &   2 &   283 & $5.573 \times 10^{16}$ \\
 6 &  1 &   1 &   294 & $3.181 \times 10^{16}$ \\
 7 &  2 &   2 &   301 & $2.628 \times 10^{16}$ \\
 8 &  1 &   1 &   305 & $1.259 \times 10^{15}$ \\
 9 &  1 &   1 &   309 & $1.772 \times 10^{16}$ \\
10 &  1 &   1 &   314 & $1.851 \times 10^{15}$ \\
\bottomrule
\end{tabular}

\bigskip
{\centering\textsc{Cumulative}\par}
\smallskip

\begin{tabular}{r rr rr}
\toprule
& \multicolumn{2}{c}{\textbf{MOEA} (pm10-L)} 
& \multicolumn{2}{c}{\textbf{Block algorithm}} \\
\cmidrule(lr){2-3} \cmidrule(lr){4-5}
\textbf{Front} $k$ 
& \textbf{Points} & \textbf{Strategies} 
& \textbf{Points} & \textbf{Strategies} \\
\midrule
 1 &  9 &   23 &    83 & $1.638 \times 10^{8\phantom{0}}$ \\
 2 & 30 &   77 &   241 & $1.396 \times 10^{13}$ \\
 3 & 43 &  107 &   466 & $1.936 \times 10^{16}$ \\
 4 & 48 &  114 &   727 & $6.650 \times 10^{16}$ \\
 5 & 50 &  116 &  1010 & $1.222 \times 10^{17}$ \\
 6 & 51 &  117 &  1304 & $1.540 \times 10^{17}$ \\
 7 & 53 &  119 &  1605 & $1.803 \times 10^{17}$ \\
 8 & 54 &  120 &  1910 & $1.816 \times 10^{17}$ \\
 9 & 55 &  121 &  2219 & $1.993 \times 10^{17}$ \\
10 & 56 &  122 &  2533 & $2.012 \times 10^{17}$ \\
\bottomrule
\end{tabular}
\caption{Distribution of results across dominance fronts for the best-performing MOEA configuration (pm10-L) and the block algorithm. For each front~$k$, we report the number of distinct objective-space points and corresponding strategies, both per front and cumulatively.}
\label{table:result_distribution}
\end{table}

\cref{table:result_distribution} confirms that most MOEA points and strategies occur on the first four fronts, with only marginal additions on later fronts. The same table also shows how these results compare with the block algorithm and reveals a significant discrepancy between the number of strategies identified by the MOEA and those found by the block algorithm. Such a discrepancy is expected: the block algorithm defines feasible strategies implicitly as Cartesian products of filtered block-level topology sets (see \cref{eq:cartasian_product}), so one objective point can represent millions of distinct strategies. The MOEA, by contrast, can only contain strategies explicitly instantiated as chromosomes, and the number present at any generation is bounded by the population size.

This difference explains why the MOEA recovers only a small fraction of the strategies associated with each objective point, even when many such strategies exist in principle. As shown in \cref{fig:front_distribution_conf010}, a substantial portion of the Pareto front produced by the block algorithm lies in the region where the number of non-reference topologies attains its maximum value. In this region, the MOEA identifies only a single solution in the fourth dominance front, indicating that this part of the search space is particularly difficult for the evolutionary approach to explore.

This behavior reflects the MOEA's difficulty in reaching extreme regions of the search space. Configurations with many non-reference blocks lie at the boundaries of the Pareto front and are therefore structurally distant in the encoding space from the strategies typically produced by crossover and mutation; as a result, they are underrepresented or missed entirely, even though the block algorithm readily identifies them. In total, 14 Pareto-optimal points correspond to configurations in which all blocks are non-reference. Excluding these boundary points leaves 69 first-front points; within this more accessible subset, the MOEA recovers 9 points, or 13\%.

\section{Discussion and conclusions}
\label{sec:conclusion}

\subsection{Summary of contributions}
This paper formalizes day-ahead transmission topology control as a sequential multi-objective optimization problem with four TSO-relevant objectives: worst-case line loading under $N-1$ security, topological depth, number of switching actions, and time spent in non-reference topologies. Two solution methods were compared on real Dutch high-voltage-grid data from TenneT.

The proposed block algorithm exploits the temporal block structure to enumerate the Pareto front exactly over the precomputed admissible topology set. With fixed operational bounds on depth and switch count, its evaluation count grows polynomially with the number of time steps, and it solves a highly congested day in under three minutes once topology-level load-flow quantities are available. 

Beyond providing optimal strategies, the block algorithm also serves as a benchmark for other methods: its Pareto front provides a reference set for distance-based metrics such as $\mathrm{IGD}^+$, and its higher-order dominance fronts let one distinguish full from partial convergence of a heuristic.

A tailored NSGA-III MOEA was developed, featuring direct variable-length encoding and structure-guided initialization across depth and switch counts. Larger populations improved front coverage, while mutation and crossover probabilities had a comparatively minor effect within the tested range. No MOEA configuration recovered solutions with $LF_1 < 100\%$, whereas the block algorithm found such solutions at low depth and with few switches. For the stated model, bounds, and case study, the block algorithm is therefore the stronger decision-support method; MOEA approaches remain relevant when additional objectives or more complex modeling assumptions break the block structure exploited here.

\subsection{Implications for grid operations} 
The results have three implications for TSO operation. First, on the studied day, strategies with $LF_1 < 100\%$ (that is, no DC thermal overloads for the considered $N-1$ contingency set) require sustained operation outside the reference topology. Second, the block algorithm yields a complete menu of Pareto-optimal plans rather than a single recommended action, allowing operators to trade switching burden against topological depth and security margin. Third, the 24-hour case is computationally tractable after topology-level load-flow preprocessing, which makes the method compatible with day-ahead congestion-management workflows.

\subsection{Limitations and future work}
Several limitations of the present study should be acknowledged, each of which suggests directions for future research.

\paragraph{Single-day evaluation}
The analysis uses one highly congested operational day, which provides a helpful stress test but does not capture variability across seasons, demand levels, and renewable generation profiles. Future work should evaluate a representative sample of operational days.

\paragraph{DC power flow approximation}
All load-flow calculations use the DC approximation~\cite{10314781}, which may diverge from full AC power flow results, potentially rendering some strategies suboptimal under more accurate models~\cite{Li2024TwoGrid}. The injections are deterministic, so forecast errors can affect the recommended strategies. Future work should incorporate AC validation and assess the impact of forecast uncertainty.

\paragraph{Preprocessing bottleneck}
The block algorithm runs quickly after preprocessing, but it requires topology-level $LF_1$ and depth values for every time step. Computing these values for all available topologies over 24 hours is substantial. 
One promising direction to mitigate this bottleneck is the development of an imitation learning agent trained on block algorithm outputs~\cite{IL_matthijs, GNN_matthijs}: by learning to map grid conditions directly to high-quality switching strategies, such an agent could enable near-real-time decision support once trained offline.

\paragraph{Objective completeness}
The four objectives are grounded in TenneT's operational practice~\cite{ViebahnGridOptions}, but they do not cover all aspects of grid performance. Notably absent is a reliability metric such as the Loss-of-Load Probability (LOLP), which quantifies the risk that demand cannot be met~\cite{Zhang2014OptimalReliability}. Integrating such metrics would yield a more comprehensive assessment of strategy quality. More fundamentally, TenneT has identified a distance-based objective that minimizes the number of busbar operations required to transition between successive topologies. Because this objective depends on pairwise relations between topologies rather than on properties of individual topologies, it cannot be precomputed per block and therefore cannot be accommodated by the current block algorithm without structural redesign. The MOEA framework remains the natural candidate for optimizing such objectives in the future.

\paragraph{Mutation operator design}
Random-reset mutation used by the MOEA directly affects the switching objective: each reset can introduce an additional topology change, making high mutation rates unattractive. At a rate of $p_m = 0.20$ applied to a chromosome of 24~genes, approximately 4.8~genes are reset on average, inflating the switching count by the same amount. Future work should test lower mutation probabilities and self-adaptive mutation control~\cite{Kramer}, eliminating the need for manual hyperparameter tuning.

\paragraph{Scalability to larger grids}
The current case study involves approximately $39{,}000$ available topologies per time step. As the number of switchable substations grows (whether through grid expansion or finer-grained substation modeling), the collection of available topologies at every time step may increase by orders of magnitude, challenging the precomputation step on which the block algorithm depends. Quantifying this scaling behavior and developing strategies to manage it (e.g., topology pre-screening heuristics) are important directions for extending the approach to larger networks.

\subsection{Concluding remarks}
Topology reconfiguration offers a low-cost, immediately deployable congestion management measure that leverages existing infrastructure. For day-ahead transmission operation, topology reconfiguration is valuable only if it can be assessed under explicit security and switching constraints. This paper shows that, for the stated DC $N-1$ thermal-security model, the combinatorial complexity of full-day topology planning can be resolved exactly by exploiting temporal block structure. The block algorithm provides TSOs with a complete set of Pareto-optimal plans over the admissible topology set and a rigorous benchmark for future heuristic or learning-based topology-control methods.

\section*{Acknowledgements}
\noindent The authors would like to thank TenneT TSO B.V. for providing access to the operational data used in this study. AI4REALNET has received funding from the European Union’s Horizon Europe Research and Innovation program under the Grant Agreement No 101119527. However, views and opinions expressed are those of the authors only and do not necessarily reflect those of the European Union. Neither the European Union nor the granting authority can be held responsible for them. The work of A.~Zocca is partially supported by the NWO Vidi grant ``\textit{Power Network Optimization in the Age of Climate Extremes}'' \href{https://doi.org/10.61686/GOOEL09973}{GrantID 10.61686/GOOEL09973}.

\section*{Data availability}
The operational data used in this study are derived from real-world grid data provided by TenneT TSO B.V. and are not publicly available. The dataset used in the analyses may be made available upon reasonable request, subject to the approval of TenneT TSO B.V.

\bibliographystyle{abbrv}
\bibliography{export.bib}

\appendix

\section{Computational complexity of the block algorithm}
\label[appsection]{app:blockalgorithmcomplexity}
The tractability of the block algorithm rests on two observations: the number of distinct ways in which the time horizon can be partitioned into blocks is relatively small (as we quantify below), and three of the four objectives take discrete integer values, which allows efficient nondominated sorting. Together, these properties enable the algorithm to pre-compute and store all relevant objective values for each block and to enumerate the complete Pareto front over the admissible topology set in polynomial time for fixed operational bounds on depth and switch count. We now make this precise by analyzing the computational complexity of the block algorithm in detail. 

A block $B$ is uniquely determined by its starting point $t_s$ and its end point $t_e$ from $T$. Therefore, the total number of unique blocks is equal to the number of ways to choose two distinct time steps from $|T| = T_{\max}$ available options, which is given by $\binom{T_{\max}}{2}$, plus the $T_{\max}$ blocks of length $1$, for which $t_s=t_e$. In our experimental setting, we will work with a one-day time horizon, subdivided into $T_{\max}=24$ unit periods, each representing an hour. For this choice, the total number of distinct blocks is quite modest and equal to $ |\mathcal{B}| = \binom{24}{2} + 24 = 300$, allowing us to pre-compute and store all objective values in \cref{sec:problem_objectives} for each block. 

As described in \cref{alg:block}, the algorithm iterates through every possible combination of $d$, $M$, and $\tildeR$, corresponding to a maximum topology depth, block configuration, and set of reference block assignments. The number of iterations is therefore based on the following components: 
\begin{itemize}
    \item The number of possible block configurations.
    \item The number of possible reference topology assignments given a block configuration.  
    \item The number of depth values considered.  
\end{itemize}
In the following paragraphs, we examine in detail how each component contributes to the total number of iterations required to find all solutions. 

The total number of block configurations depends on the number of possible switching time steps.
Assuming we start with a block at time step $t=0$, there are $T_{\max} - 1$ potential switching time steps within the time horizon. At each time step, there are two possibilities: start a new block or ``extend'' the current block. As a result, the total number of block configurations is equal to $2^{T_{\max}-1}$. Many of the resulting block configurations prescribe a large number of switching actions, possibly even at every time step. In practice, however, we are interested in configurations with a limited number of switches, in line with objective 3 in \cref{sec:problem_objectives}.

For any $l \in \{0, \dots, T_{\max}-1\}$, let $p(l)$ be the number of block configurations with exactly $l$ switching time steps and thus consisting of exactly $l+1$ blocks, which is equal to
\begin{equation}
\label{eq:switch_exact}
    p(l) \coloneq \binom{T_{\max} - 1}{l}.
\end{equation}
Consequently, the number of block configurations with at most $l$ switches is given by
\begin{equation}
\label{eq:switch_cumulative}
    P(l) \coloneq \sum_{i = 0}^{l} p(i) = \sum_{i = 0}^{l} \binom{T_{\max} - 1}{i}.
\end{equation}
Note that calculating $P(T_{\max}-1)$, we recover $2^{T_{\max} - 1}$, which is the total number of block configurations.

Next, we consider the number of possible reference-topology assignments for a given block configuration. In each block $B \in M$, we can either use the reference topology or not. Thus, the total number of ways to assign the reference topology blocks to a configuration $M$ is equal to $2^{|M|}$. However, the number of possibilities to evaluate can be reduced by considering only combinations in which no two consecutive (connected) blocks are selected as reference blocks, as explained below.

We say that two consecutive blocks, $B_i$ and $B_{i+1} \in M$, are \textit{connected} if the end point of $B_i$ directly precedes the start point of $B_{i+1}$, that is, if $t_{e_i} + 1 = t_{s_{i+1}}$. Suppose both connected blocks are selected as reference blocks. In that case, the block configuration $M$ is equivalent to the configuration $M'$ obtained by replacing $B_i$ and $B_{i+1}$ with a single reference block $B^{*} = B_i \cup B_{i+1} = \{t_{s_i},\dots, t_{e_{i+1}}\}$. Selecting $B_i$ and $B_{i+1}$ together or selecting the merged block $B^{*}$ yields the same total number of time steps outside the reference topology,
$
z(\tildeR) = z(\tildeR'),
$
and the same $LF_1$ value,
$
LF_1(d, M, \tildeR) = LF_1(d, M', \tildeR'),
$
where $\tildeR' = (\tildeR \setminus \{B_i, B_{i+1}\}) \cup \{B^{*}\}$.  

The block algorithm processes configurations in increasing order of $|M|$; that is, it first evaluates all configurations consisting of a single block, then all configurations with two blocks, and so on.  
Since $M'$ contains one fewer block than $M$, i.e., $|M'| = |M| - 1$, this reduced configuration corresponds to a case that has already been evaluated in an earlier iteration of the algorithm and therefore does not need to be reevaluated.

For a block configuration $M$ with $N=|M|$, the feasible reference assignments $\mathcal{A}(M)$ in~\cref{eq:block2} are exactly the subsets of $\{1,\dots,N\}$ with no adjacent indices. If $f(N)$ denotes their number, conditioning on whether block $N$ is selected gives
$
f(N)=f(N-1)+f(N-2),
$
with $f(0)=1$ and $f(1)=2$. Hence $f(N)=F_{N+2}$, where $F_0=0$, $F_1=1$, and $F_{k+2}=F_{k+1}+F_k$ denotes the Fibonacci sequence. Since a configuration with $l$ switching time steps has $l+1$ blocks, it admits $F_{l+3}$ feasible reference assignments.

The final component of the computational complexity concerns depth.
For each block configuration $M$, the algorithm iteratively filters the topology sets $G_B$ for all $B \in M$ across depth levels $d = 1, 2, \dots, d^{\max}$.  
At each depth level, only the filtered subsets $G_B^{(\leq d)}$ are retained, and the corresponding $LF_1$ for the block configuration and each block reference assignment needs to be reevaluated.

Knowing the computational complexity of each component, we can express the total number of iterations required for a given time horizon, maximum depth, and maximum number of switches (equivalently, a maximum number of blocks minus one) using the following equation: 
\begin{equation}
\label{eq:switch_cumulative_fibonacci}
    N_{\text{eval}}(d^{\max}, l^{\max}, T_{\max}) 
    = (d^{\max}+1) \cdot 
    \sum_{l = 0}^{l^{\max}} 
    \binom{T_{\max} - 1}{l} \, F_{l+3},
\end{equation}
where $d^{\max}$ and $l^{\max}$ denote the maximum depth and switches, respectively.  

The total number of evaluations grows linearly with the number of depth levels and polynomially with $T_{\max}$ of degree $l^{\max}$, leading to an overall complexity of $\mathcal{O}\!\left((d^{\max}+1)\,T_{\max}^{\,l^{\max}}\right)$.
In practice, $d^{\max}$ and $l^{\max}$ are small constants determined by operational constraints, meaning the algorithm scales polynomially with the number of time steps for fixed operational bounds on depth and switch count and remains computationally tractable for realistic problem sizes. If $l^{\max}$ were allowed to increase with $T_{\max}$, the complexity would become exponential, but such configurations are not relevant for operational use cases.

\section{MOEA implementation and evaluation details}
\label[appsection]{app:moea}

\subsection{Representation}
\label[appsection]{app:chromosome_formulation}
In an evolutionary algorithm, each candidate solution is encoded as a chromosome (genotype); the mapping to the actual decision variables (phenotype) determines how genetic operators translate into meaningful changes~\cite{EibenNaturalComputing}. Our choice of chromosome formulation is guided by~\cite{Beke2021AMultigraphs}, which evaluates three representation schemes for multi-objective sequential decision problems: direct variable-length encoding, integer-valued priority-based encoding, and random-key-based encoding. Their empirical comparison shows that random-key-based encoding yields the strongest overall performance across multigraph shortest-path instances, owing to its flexibility and favorable search dynamics. However, for the congestion problem introduced in \cref{sec:problem_objectives}, this approach becomes impractical. In our setting, each gene must encode a topology–time pair, implying a chromosome length of
$
    \sum_{t \in T} { |G_t|},
$
which can quickly get very large in our setting. Such long chromosomes incur high memory costs and increased runtime in both variation and evaluation. The same issue arises with the integer-valued priority-based representation, whose length scales in the same way, rendering it similarly unsuitable for our problem structure. Direct variable-length encoding is therefore the natural choice. In this representation, the number of genes scales directly with the number of time steps, resulting in compact chromosomes and efficient operator application. An additional advantage is the one-to-one correspondence between genotype and phenotype: each chromosome encodes a unique solution trajectory. 
The authors in~\cite{Beke2021AMultigraphs} highlight that this bijectivity avoids representational redundancy: in many-to-one encodings, multiple genotypes can decode to the same phenotype and therefore share identical fitness values, reducing search efficiency.

A commonly cited drawback of direct variable-length encoding is the potential formation of loops when applying crossover or mutation, which may require explicit repair procedures to restore feasibility~\cite{Beke2021AMultigraphs}. 
In our setting, however, loop formation is structurally impossible. A solution consists of a topology sequence $s = (g_t)_{t \in T}$, where each topology corresponds to a distinct time step. Transitions are allowed only from one time step to the next and never to a previous or current time step. Therefore, the problem's structure prevents loops by design. For these reasons, direct variable-length encoding has been adopted as the chromosome representation.

\subsection{Initial Population}
\label[appsection]{sec:initial_population}
Before the evolutionary process begins, the initial population must be generated. There are two main prerequisites for this population: each chromosome must encode a feasible strategy, and the solutions should be sufficiently diverse within the objective space. Diversity in the initial population improves coverage of the search landscape, enabling evolutionary operators to explore multiple promising regions and reducing the risk of premature convergence \cite{EibenNaturalComputing}.

A straightforward approach is random initialization, which in our setting amounts to independently selecting a topology $g_t \in G_t$ uniformly at random for each time step $t \in T $. While simple, this method performs poorly for our problem for two structural reasons. First, the topological depth distribution is intrinsically skewed. Indeed, by construction, the set $G_t$ of available topologies at time step $t$ contains substantially more topologies with a large topological depth. Uniform random sampling is therefore unlikely to generate strategies with small depth values, which are preferable since one of the main objectives (see \cref{sec:problem_objectives}) is to minimize the maximum topological depth along a strategy. Second, on average, the samples have an excessive number of switches. Indeed, combining topologies independently sampled for each time step $t$ tends to produce strategies with many switches, far more than would occur in realistic or high-quality solutions.

As a consequence, a uniform random initialization yields a population that is diverse in its raw encodings but not diverse in the relevant objective-space dimensions (number of switches and depth). This hinders the chosen evolutionary algorithm's ability to produce a high-quality approximation of the Pareto front. To address this limitation, we design a custom sampling procedure that systematically balances diversity across these structural properties of strategies.

\underline{Structure-Guided Sampling Based on Switch Count}. For any given switch count $l$, we first select a random subset of block configurations $\mathcal{M}'(l) \subset \mathcal{M}(l)$, defined earlier in~\cref{eq:block_conf_sets}. For each block configuration $M \in \mathcal{M}'(l)$, we consider its blocks $B \in M$ and for each block $B$ we draw one topology at random from the corresponding topology collection $G_B$. This procedure generates strategies that share the same number of switches $l$ but differ in their switching locations and topology choices.
Formally, the set of all possible strategies with exactly $l$ switches is
\[
\mathcal{S}^{(l)} := \bigcup_{M \in \mathcal{M}(l)} \left\{ s = (g_B)_{B \in M} ~:~ g_B \in G_B,\; \forall B \in M \right\},
\]
from which we sample a finite subset $\mathcal{P}^{(l)} \subset \mathcal{S}^{(l)}$. The switch-based component of the initial population is then obtained by merging all these subsets into
\[
\mathcal{P}_l =
\bigcup_{l = 0}^{l^{\max}} \mathcal{P}^{(l)}.
\]
This ensures that the initial population contains strategies spanning a broad yet well-controlled number of switches.

\underline{Structure-Guided Sampling Based on Depth}. To introduce diversity for the maximum depth objective, we similarly generate strategies associated with controlled maximum depth values. The depth-based component of the initial population is defined as
\[
\mathcal{P}_d =
\bigcup_{d = 0}^{d^{\max}} \mathcal{P}^{(d)},
\]
where $\mathcal{P}^{(d)}$ is a random subset of feasible strategies whose maximum depth is exactly equal to $d$. Because the depth of a strategy is determined after sampling, a small number of sampled strategies may achieve depth values lower than the target depth $d$. These outliers are removed and resampled to maintain the intended distribution across depth levels. The final initial population combines diversity across both switch count and maximum depth via
\[
\mathcal{P} = \mathcal{P}_l \cup \mathcal{P}_d.
\]
\underline{Structure-Guided Sampling Sizes}. The sizes of the individual subsets used to construct the initial population are controlled by two predefined sampling parameters, denoted by $\bar{l}$ and $\bar{d}$.  
Specifically, for the switch-based sampling, each subset $\mathcal{P}^{(l)}$ associated with a switch count $l$ is sampled with fixed cardinality
\[
\lvert \mathcal{P}^{(l)} \rvert = \bar{l}, \qquad \forall\, l \in \{0,\dots,T_{\max}-1\}.
\]
Similarly, for the depth-based sampling, each subset $\mathcal{P}^{(d)}$ associated with a target maximum depth $d$ is sampled with fixed cardinality
\[
\lvert \mathcal{P}^{(d)} \rvert = \bar{d},
\qquad \forall\, d \in \{1,\dots,d^{\max}\}.
\]
The case $d = 0$ is excluded from this sampling procedure, since there is a unique topology at that depth, i.e., the reference topology $\tilde{g}$, which is always added to the initial population.

Parameters $\bar{l}$ and $\bar{d}$ directly determine the contribution of the switch-based and depth-based components to the total population size $\mathcal{P}$, yielding
\begin{equation}
\label{eq:initial_population_size}
    |\mathcal{P}| =  T_{\max} \cdot \bar{l} + d^{\max} \cdot \bar{d} +1.
\end{equation}

This construction yields an initial population that is (i) composed exclusively of feasible strategies and (ii) explicitly diverse in the objective-space dimensions that govern search difficulty. As a result, it provides a more informative and balanced starting point for the MOEA than uniform random sampling, thereby improving its exploratory capabilities and overall optimization performance.

\subsection{Variation: Crossover and Mutation}
\label[appsection]{app:moea_variation}
This subsection describes the variation operators, crossover, and mutation, used in the proposed MOEA algorithm. Variation operators are central to evolutionary algorithms: crossover primarily supports exploitation by recombining favorable genetic material from parent solutions, whereas mutation promotes exploration by introducing stochastic perturbations that enable the discovery of novel regions of the search space.
A key design requirement for the variation operators is that the resulting offspring remain feasible. Recall that at each time step $t$, only a restricted collection $G_t$ of topologies is available. A chromosome encodes a strategy $s = (g_t)_{t \in T}$ subject to the constraint that each gene satisfies $g_t \in G_t$ for all $t \in T$. Any variation operator must therefore ensure that this constraint is satisfied by all offspring.

A crossover method that ensures this feasibility is the \textit{$k$-point crossover operator}. In this operator, $k$ crossover points are selected uniformly at random along the chromosome, partitioning it into $k+1$ contiguous segments. The offspring's chromosomes are then formed by alternately copying segments from the two parents. Crucially, because each gene position $t$ in the offspring is copied from the corresponding position in one of the parents, and both parents satisfy $g_t \in G_t$, feasibility is preserved by construction. 

Beyond feasibility preservation, the choice of $k$-point crossover is further motivated by the temporal structure of the decision horizon. For example, in uniform crossover, where the allele at each time step is selected independently from either parent, the recombination process introduces many topology changes in the offspring. As a result, contiguous temporal patterns present in the parents become frequently fragmented. In contrast, $k$-point crossover preserves contiguous temporal segments of the chromosome, thereby maintaining the integrity of temporal sequences inherited from the parents. This results in fewer disruptions of temporal structure, so offspring tend to retain effective switching patterns. In conclusion, $k$-point crossover provides a simple and widely used recombination mechanism that respects temporal coherence while preserving feasibility, without introducing additional problem-specific assumptions or implementation complexity.

Mutation is implemented via a \textit{random-reset operator}. Each gene $g_t$ is independently selected for mutation with probability $p_m$. If a mutation occurs at time step $t$, the current gene is replaced by a new topology sampled uniformly at random from the set $G_t$. Since all possible replacements are in $G_t$, this operator also ensures feasibility. 

\subsection{Evaluation and Selection}
\label[appsection]{app:moea_selection}
The evaluation and selection mechanisms jointly drive the evolutionary search toward high-quality solutions. Evaluation assigns fitness values to chromosomes based on their performance with respect to the problem objectives, while selection determines which individuals are retained and propagated to subsequent generations.
The objective functions defined in~\cref{sec:problem_objectives} are used directly as fitness functions, which is possible because the genotype and phenotype coincide for the chosen direct encoding (see~\cref{app:chromosome_formulation}). Selection must balance convergence toward high-quality solutions with diversity preservation across four competing objectives.

Since this qualifies as a many-objective problem, we adopt the NSGA-III framework~\cite{Deb2014AnConstraints}, which extends the classical NSGA-II framework by combining nondominated sorting with reference-point-based diversity preservation. 
In nondominated sorting, solutions are ranked into successive fronts based on Pareto dominance: the first front contains all solutions that are not dominated by any other solution, the second front contains solutions dominated only by those in the first front, and so on.
This approach does not require any weighting schemes and thus avoids objective aggregation, treating all objectives equally. Moreover, the explicit Pareto-ranking structure of NSGA-III facilitates the interpretability of the trade-offs among competing objectives. From a computational perspective, NSGA-III has the same asymptotic per-generation complexity as NSGA-II, with only minor additional overhead due to reference-point association.
NSGA-III has been shown to outperform MOEA/D~\cite{Chen2021MOEADOptimization}, IBEA~\cite{Zitzler2004IndicatorBasedSearch}, and MOPSO~\cite{Coello2004HandlingOptimization,Coello2002MOPSO} on combinatorial many-objective benchmarks~\cite{Ma2023AApplications}. In view of all these considerations, NSGA-III is adopted as the selection algorithm in this study.

\subsection{Hyperparameters for MOEA}
\label[appsection]{app:hyperparameters}
The proposed MOEA is inherently stochastic, so each experimental configuration is evaluated over 15 independent random seeds. Since the runs are independent, they are executed in parallel without increasing wall-clock time.

To assess the algorithm's robustness and sensitivity to key design choices, several hyperparameter configurations are considered. The experimental variations focus on three central parameters: population size, crossover probability $p_c$, and mutation probability $p_m$. In all configurations, the crossover and mutation probabilities are set as complements, $p_c = 1 - p_m$, so that the experimental variations effectively reduce to a single parameter. Note that $p_c$ governs whether crossover is applied to a 
selected pair of parents, while $p_m$ is the per-gene probability of random reset; the two operators act at different levels and their complementary coupling is a design choice rather than a theoretical requirement.
All remaining algorithmic components are kept fixed across experiments to isolate the effects of these parameters. Evaluating multiple hyperparameter combinations across several random seeds increases the likelihood of obtaining a well-distributed and high-quality approximation of the Pareto front. \cref{tab:experiment_parameters} summarizes the algorithmic settings that are fixed in all experiments. 

\begin{table}[h!]
\centering
\begin{tabular}{l|l}
\textbf{Parameter} & \textbf{Value/choice} \\
\hline \hline
Seeds & 15 \\
Crossover method & $k$-point crossover, with $k=2$ \\
Mutation method & Random reset \\
Selection method & NSGA-III \\
Reference directions & 100 \\
Termination & 30 minutes \\
\end{tabular}
\caption{Standard algorithm configuration}
\label{tab:experiment_parameters}
\end{table}

One requirement for NSGA-III is the reference directions. The reference directions are uniformly distributed across the normalized action space and are generated using a Riesz s-energy minimization approach \cite{HARDIN2005174}. For the four-objective problem, $100$ reference directions are used, as this number provides sufficiently fine coverage of the normalized objective space while remaining compatible with the chosen population sizes.

All experiments are run on a Microsoft Dev Box virtual machine equipped with 32 GB of RAM and a 64-core processor (2.44\,GHz). The block algorithm described in \cref{sec:block_algo} computes the full Pareto front for a problem instance of the considered size in under three minutes. To ensure that the evolutionary algorithm is not constrained by limited runtime and has ample opportunity to explore the search space, each MOEA run is allocated a substantially larger computational budget of 30 minutes. Note that, given a fixed time budget, mutation and crossover probabilities have no significant effect on the number of generations; only the population size does.

The first configuration-specific hyperparameter is the population size, which strongly influences the balance between exploration and convergence. Smaller populations allow more generations to be completed within a fixed time budget, increasing the frequency of crossover and mutation events, and enabling faster short-term adaptation. However, they may suffer from limited diversity and a higher risk of premature convergence. Larger populations preserve greater genetic diversity and improve coverage of the Pareto front, but typically evolve more slowly because fewer generations can be completed. To examine this trade-off, experiments are conducted using small (\textbf{S}), medium (\textbf{M}), and large (\textbf{L}) population sizes. 

Each population size is determined by the sampling parameters $\bar{l}$ and $\bar{d}$ introduced in \cref{sec:initial_population}. For configuration \textbf{S}, the sampling parameters $\bar{l}$ and $\bar{d}$ are both set to 30, yielding an initial population of 811 individuals. We note that, in this way, the population size is set to be an order of magnitude higher than the number of solution points found by the block algorithm (83 in total, as detailed later in \cref{table:result_distribution}), hence giving the MOEA a fair chance to discover all these optimal solutions. The parameters $\bar{l}$ and $\bar{d}$ are linearly scaled and set to 45 for configuration \textbf{M} and 60 for configuration \textbf{L}. This results in total population sizes of 1,216 and 1,621, respectively, as summarized in \cref{tab:number_of_generations}.

For each hyperparameter configuration, we run 15 independent evolutionary runs in parallel, each with a different seed. Although the termination criterion is time-based in principle, in order to get comparable results across different seeds, individual runs are, in practice, terminated after a fixed number of generations, as detailed below in~\cref{tab:number_of_generations}.

\begin{table}[!ht]
\centering
\begin{tabular}{c|c|c}
\textbf{Configuration} & \textbf{Population size} & \textbf{Target \# generations} \\
\hline \hline
\textbf{S} & 811 & 9230  \\
\hline
\textbf{M} & 1216 & 6340 \\
\hline
\textbf{L} & 1621 & 4700 \\
\end{tabular}
\caption{Population size and target number of generations for \textbf{S}, \textbf{M}, and \textbf{L} configurations.}
\label{tab:number_of_generations}
\end{table}

In addition to population size, the mutation and crossover probabilities directly influence the balance between exploration and exploitation in the evolutionary process. Mutation introduces stochastic perturbations that promote exploration of previously unvisited regions of the search space, while crossover exploits information from existing solutions to accelerate convergence toward high-quality regions. The balance between these operators determines whether the algorithm emphasizes exploration or exploitation, thereby affecting convergence speed, solution diversity, and robustness to local optima. To study these effects, several combinations of $p_m$ and $p_c$ are evaluated for each population size. The resulting configurations are listed in \cref{tab:algorithm_configurations}.

\begin{table}[!ht]
\centering
\begin{tabular}{l|c|c|c}
\textbf{Configuration} & \textbf{Population size} & \textbf{Mutation prob.} $p_m$ & \textbf{Crossover prob.} $p_c$ \\
\hline \hline
pm05-S & 811  & 0.05 & 0.95 \\
pm05-M & 1216 & 0.05 & 0.95 \\
pm05-L & 1621 & 0.05 & 0.95 \\
\hline
pm10-S & 811  & 0.10 & 0.90 \\
pm10-M & 1216 & 0.10 & 0.90 \\
pm10-L & 1621 & 0.10 & 0.90 \\
\hline
pm15-S & 811  & 0.15 & 0.85 \\
pm15-M & 1216 & 0.15 & 0.85 \\
pm15-L & 1621 & 0.15 & 0.85 \\
\hline
pm20-S & 811  & 0.20 & 0.80 \\
pm20-M & 1216 & 0.20 & 0.80 \\
pm20-L & 1621 & 0.20 & 0.80 \\
\end{tabular}
\caption{Algorithm configurations used in the experiments. The identifier encodes the mutation probability $p_m$ ($\texttt{pmXX} = p_m = \texttt{XX}/100$) and the population sizes corresponding to configurations \textbf{S}, \textbf{M}, and \textbf{L}.}
\label{tab:algorithm_configurations}
\end{table}

For each algorithmic configuration, populations obtained from different random seeds are combined prior to performance evaluation. To formalize this aggregation, we introduce the following sets: $E$, the finite set of random seeds, and $R$, the finite set of generations. We denote by $P_{r,e}$ the population at generation $r \in R$ for seed $e \in E$, where each element $s \in P_{r,e}$ represents a single strategy. The final generation index is denoted by $\tilde{r}$. The combined population at generation $r$ is then defined as
\[
\overline{P}_{r} \coloneq \bigcup_{e \in E} P_{r,e}.
\]
Although all weakly nondominated solutions are retained during the evolutionary process, Pareto front analysis is performed only on unique points. As a result, the reported fronts effectively consist of strongly nondominated strategies, ensuring that duplicate or equivalent solutions do not bias the evaluation. The Pareto front at generation $r$ is then given by $PF(\overline{P}_r)$, which represents the set of strategies that are nondominated with respect to all strategies generated across different random seeds at that generation.

\subsection{Evaluation Measures for MOEA}
\label[appsection]{app:evaluation measures}
To evaluate the MOEA's performance throughout the evolutionary process, two quantitative measures are employed. The first measure assesses convergence and diversity with respect to the optimal solution set, whereas the second measure evaluates how well the algorithm reproduces different dominance levels identified by the reference method, i.e., the block algorithm.

The first metric is the \textit{Inverted Generational Distance Plus} ($\mathrm{IGD}^+$), which measures how closely the obtained solutions approximate the true Pareto front, i.e., the set of solutions generated by the block algorithm. Formally, $\mathrm{IGD}^+$ is computed by averaging, over all reference points on the true front, the minimum Euclidean distance to the set of solutions obtained by the MOEA, using a dominance-preserving distance measure \cite{Ishibuchi2019ComparisonSolutions}. Lower values of $\mathrm{IGD}^+$ indicate better convergence toward and coverage of the true Pareto front.

Before computing $\mathrm{IGD}^+$, all objective values are normalized using min--max normalization. The normalization is performed with respect to the lower and upper values used to scale each objective. For the discrete objectives, that is, topological depth, number of switches, and number of non-reference topologies, these values are determined by the operational limits used to construct the exact reference fronts. For the continuous objective $LF_1$, the values are obtained from the best and worst feasible $LF_1$ values produced by the block algorithm given the whole objective space (before nondominated sorting). The resulting normalization values are summarized in \cref{nadir}.

\begin{table}[!ht]
\centering
\begin{tabular}{c|c|c|c|c}
\textbf{Objectives} & \textbf{Depth} & \textbf{Switches} & \textbf{Non-ref.~steps} & \textbf{$LF_1$} \\
\hline \hline
\textbf{Ideal} & 0 & 0 & 0 & \, 77.3 \\
\hline
\textbf{Maximum} & 3 & 5 & 24 & 212.5 \\
\end{tabular}
\caption{Ideal and maximum objective values defining the bounds of the objective space.}
\label{nadir}
\end{table}

The maximum values for depth and switches directly reflect the operational bounds imposed on the exact reference calculation: $d^{\max} = 3$ corresponds to the maximum topological depth available in the dataset (see \cref{tab:values_counts}), and $s^{\max} = 5$ is the maximum number of switching time steps represented in the reference fronts. These same bounds apply to the MOEA, ensuring that both methods operate over an identical search space.

Beyond the Pareto front itself, the block algorithm also identifies solutions belonging to subsequent dominance levels. The first level corresponds to the Pareto front, while the second level contains solutions dominated only by those in the first level, and so forth. Although solutions in higher dominance levels are suboptimal, they may still represent competitive trade-offs that are close to optimal with respect to specific objectives. Incorporating these additional fronts enables a more nuanced assessment of MOEA performance: even when the algorithm fails to recover all Pareto-optimal solutions, identifying solutions in lower dominance levels indicates partial convergence and meaningful progress toward optimality.

To quantify how well the MOEA reproduces these reference dominance fronts, we define the following performance measures:
\begin{equation}
\label{eq:evaluation_metrics}
I_k(P) := |\, PF(P) \cap F_k \,| \quad \text{ and } \quad \hat{I}_k(P) := \frac{|\, PF(P) \cap F_k \,|}{|\, F_k \,|},
\end{equation}
where $F_k$ denotes the $k$-th reference front obtained from the block algorithm, ordered by increasing dominance rank. 
The quantity $I_k(P)$ counts the number of solutions in the MOEA Pareto front $PF(P)$ that belong to the reference front $F_k$, while $\hat{I}_k(P)$ normalizes this count by the size of $F_k$, yielding a relative coverage measure. Together, these metrics provide insight into both the absolute and relative ability of the MOEA to find solutions at different dominance levels. In this study, we consider $k=10$ different dominance levels. 

Note that the cardinality of $I_k(P) $ depends on how equality between solutions is defined across objectives. In particular, $LF_1$ is the only continuous objective, whereas all other objectives are discrete. As a result, exact equality in $LF_1$ values between MOEA solutions and reference solutions is unlikely, and without any numerical tolerance or rounding, the intersection $PF(P) \cap F_k $ may be empty or contain only very few solutions, even when solutions are effectively equivalent in practice. To address this issue, $LF_1$ values are first rounded. In our study, we round the $LF_1$ values to the first decimal place.

Both $\mathrm{IGD}^+$ and the dominance-level coverage measures in~\cref{eq:evaluation_metrics} are computed per generation using the combined Pareto front $PF(\overline{P}_r)$. This allows the algorithm's convergence behavior to be tracked over time, revealing not only the final solution quality but also the convergence rate and stability throughout the evolutionary process.
\end{document}